\documentclass[14pt,a4paper]{article}
\usepackage{version}
\usepackage{float}
\usepackage{amsmath}
\usepackage{amstext}
\usepackage{tabularx}
\usepackage[integrals]{wasysym}
\usepackage{marvosym}
\usepackage{manfnt}
\usepackage{ifsym}
\usepackage[psamsfonts]{amssymb}
\usepackage[unicode,bookmarks,bookmarksopen,bookmarksopenlevel=0,colorlinks,%
plainpages=false,pdfpagelabels,hypertexnames=false,pageanchor=true,breaklinks=true,%
linkcolor=blue,citecolor=blue,urlcolor=red]{hyperref}
\usepackage[left=2.5cm, right=2.5cm, top=2.5cm, bottom=2.5cm,
includehead]{geometry}
\usepackage{upref}  
\usepackage{indentfirst}
\newif\ifpdf
\ifx\pdfoutput\undefined
\pdffalse 
\else
\pdfoutput=1 
\pdftrue \fi
\ifpdf
\usepackage[pdftex]{graphicx}
\usepackage{graphicx}
\else
\usepackage{graphicx}
\fi

\newcommand\pr{\partial}
\newcommand\ben{\begin{enumerate}}
\newcommand\een{\end{enumerate}}
\newcommand\bed{\begin{itemize}}
\newcommand\eed{\end{itemize}}
\newcommand\bei{\begin{description}}
\newcommand\eei{\end{description}}
\newcommand\beear{\begin{eqnarray}}
\newcommand\eear{\end{eqnarray}}
\newcommand\beq{\begin{eqnarray*}&}
\newcommand\eeq{&\end{eqnarray*}}
\newcommand\pht{\hat{p}}

\newcommand\bpx{\overline{p}^\xi}
\newcommand\pch{\check{p}}

\newcommand\efr{\mathfrak{e}}
\newtheorem{prpsn}{Propopsition}
\newtheorem{lemma}{Lemma }
\newtheorem{theorem}{Theorem}
\newtheorem{dfntn}{Definition}

\newcommand\const{\operatorname{const}}
\newcommand\Div{\operatorname{div}}

\newcommand{\byd}{\stackrel{\mathrm{def}}{=}}

\newcommand\va{\boldsymbol{a}}
\newcommand\vb{\boldsymbol{b}}
\newcommand\vc{\boldsymbol{c}}
\newcommand\vf{\boldsymbol{f}}

\newcommand\vtheta{\boldsymbol{\theta}}
\newcommand\vu{\boldsymbol{u}}

\newcommand\vq{\boldsymbol{q}}

\newcommand\ve{\boldsymbol{e}}

\usepackage{fancyhdr}
\usepackage{lastpage}
\begin{document}
\title{ {\small \bf Prey-taxis VS a Shortwave External Signal in  Multiple Dimensions}}
\author{{\small Andrey Morgulis\footnote{ORCID 0000-0001-8575-4917X, abmorgulis@sfedu.ru}}, \\
{\footnotesize I.I.Vorovich Institute for Mathematic, Mechanics and Computer Science,}\\
{\footnotesize Southern Federal University, Rostov-na-Donu, Russia;\quad}\\
{\footnotesize Southern Mathematical Institute of VSC RAS, Vladikavkaz, Russia}
 \and 
  {\small Karrar H. Malal}\\
 {\footnotesize I.I.Vorovich Institute for Mathematics, Mechanics and Computer Sciences,}\\
{\footnotesize Southern Federal University, Rostov-na-Donu, Russia} \\
  }
\date{{\footnotesize August 13, 2025}}
\maketitle
\abstract{We consider a model of the predator--prey community with prey-taxis.  By that we mean   the capability of the predators to get moving in a certain direction on the macroscopic level in response  to the prey density gradients. Additionally, we suppose the same kind of sensitivity with respect to  one more signal, called external, the production of which goes on independently of the community state. Such a signal can be due to the spatiotemporal inhomogeneity of the environment that results from the natural or artificial reasons.  The model employs the Patlak--Keller--Segel law  for responses to both ones.  We assume that the  external signal takes a general short-wave form, and  we construct the complete  asymptotic expansions of the    short-wave solutions.  This result generalizes the prior one by Morgulis \& Malal (2025)  in two respects. First, we have addressed the case of multiple dimensions. Second, we have got rid of  assuming  the signal and corresponding solutions to take  the form of a traveling  wave,   that makes our result novel even in one dimension. Further, we apply  the short wave asymptotic to studying the stability or instability imposed by the external signal  following Kapitza' theory for upside-down pendulum. 
}



\section{Introduction}\label{ScIntr}
This article continues the studies published in articles \cite{AM1},\cite{AM2} and \cite{AM3}, and, in effect, it is  the second part of the  last one,  which include an extended introductory part. So we allude to it for guiding through the topic and the references, while putting quite a sketched  introduction here. 

All the studies mentioned (including the present one) address   the predator-prey models with the so-called prey-taxis, by which we mean the ability of the predators to move deterministically  for getting the prey. However, the models of prey-taxis dealt with there differ one to another one. At the same time,   there is a common  goal for them to pursue:  to investigate  the response to  shortwave forcing using the asymptotic techniques (homogenization, multiple scale approach, etc, see \cite{Allr-1}, for example), including the imposed stability or instability of the homogeneous equilibria in spirit of the Kapitza's theory for the upside-down  pendulum, \cite{LndLfs}.  This is what differs  the mentioned studies from many others that also address  the instabilities and the pattern formation, but with no external forcing, e.g. \cite{BrzKrv}-\cite{Grchv}

The current study rely on the classical Patlak-Keller-Siegel (PKS) formulation for the tactical flux. However, we consider the flux driven by coupling of two signals. First one is  the prey density, and the second one, called external,  is independent of the community state, and specifying its intensity is a part of the model setup. For instance, this can be  the temperature distribution over some  biotop,  or an intensity of  physical or chemical  action  exerted  from a researcher to some    community \emph{in vitro}  while doing some experiment. It is quite logical to allow sometimes for the very different scales of the signals. For instance, aside from the mentioned studies, articles \cite{Chdh1}-\cite{Chdh3} address two-scale models of ecological systems including the predator-prey one.

So, we assume the external signal to have a shortwave shape and then proceed in a way outlined above, but with no assumptions to which  study \cite{AM3} relies essentially. We neither set the spatial dimension to one nor set the signal or the solution in question to be a short traveling wave,  and, eventually,  achieve  (i) developing full asymptotic expansion for the solutions to the predator-prey system with Keller-Siegel's responding to the prey density and to a shortwave external signal in multiple dimensions with no special assumptions on the shapes of external signal or solution; (ii) finding the criteria for the imposed stability or instability due to a general shortwave external signal in multiple spatial dimensions.

The article is organized as follows. In Section~\ref{ScSttng}, we formally set the model. In Section~\ref{ScPrlmnr} we 
introduce some auxiliary matters. In Section~\ref{ScAsmpt}, we derive the asymptotic and discuss it's relation  to the particular case addressed in \cite{AM3}. In Section~\ref{ScStb}, we discuss the obtained result and its outcomes regarding the imposed stability.  
\section{Formulation of the problem}
\label{ScSttng}
\noindent 
We address  a multiple-dimensional model  that read as
\begin{eqnarray}
 &p_t+\Div{\vq}=pf(p,s),&
 \label{EqPrdTrns}\\
 &{\vq}=p\nabla\left(\chi s +\varkappa h -\mu \ln  p\right),&
  \label{EqFlx}\\
  &s_t=\delta\Delta s+sg(p,s). &
\label{EqPryDnst}  
\end{eqnarray}
Here {$t$} is time, {$x\in\mathbb{R}^n$} is  spatial coordinate, $p=p(x,t)$, ${ \vq}={\vq}(x,t)$ -- density  and   flux of the predators, and function $s=s(x,t)$ denotes the prey density. Functions $f$ and $g$ are to determine the system kinetics. We assume  that they are known  and analytic in the cone $\{p>0,s>0\}$.  Inside this cone equation \eqref{EqFlx} determines  nothing else than  PKS-flux due to the couple of signals one of which   is the prey density, $s$, and the other one is produced  externally with a known intensity denoted as $h$. 
Coefficients $\chi$ and $\varkappa$ are to measure  the predators sensitivity to the prey density  and to the external signal, correspondingly, and coefficients $\mu$ and $\delta$ are diffusivity coefficients for the predators and prey, correspondingly. They all are assumed to be constant. 
In this way, equations \eqref{EqPrdTrns}, \eqref{EqFlx} and \eqref{EqPryDnst} constitute the PKS model for the prey-taxis with external signal, or, briefly, PKS-model. We regard this system as a dimensionless one, see \cite{AM3} for a discussion on  scaling.


In the just mentioned article  Morgulis \& Malal have addressed the short traveling waves in one dimension driven by an external signal, which itself is a short traveling wave  -- that is,
\begin{eqnarray}&\label{ShrtWvSgnl0}
  h=h(x,t,\eta),\quad p=p(x,t,\eta),\quad q=q(x,t,\eta),\quad s=s(x,t,\eta), \quad x\in \mathbb{R},\ \eta=(x-ct)/\delta, \ c=\const,\quad  
&\end{eqnarray}
where $\delta$ is small parameter, $h$ is smooth in the triad of  variables $x,t,\eta$, and 
\beq 
h(x,t,\eta+2\pi)=h(x,t,\eta) \ \forall\, x,t,\eta,\quad \int\limits_0^{2\pi} h(x,t,\eta)=0\ \forall x,t.
\eeq
Here  we consider more general  signal. We set
\begin{eqnarray}&\label{ShrtWvSgnl}
h=h(x,t,x/\delta,t/\delta), 
&\end{eqnarray}
where function $h=h(x,t,\xi,\tau)$ is $\mathrm{C}^{\infty}$-smooth on $\mathbb{R}^{n+1}\times\mathbb{R}^{n+1}$ and   periodic in coordinates $\xi_1,\ldots\xi_n,\tau$ with periods  $\ell_1,\ldots\ell_n,\ell_0$, and  positive number  $\delta$ is small.    Additionally,  we assume that the average value of function $h$ over periods' box is equal to zero.

We'll be  seeking for the power asymptotic that reads 
\begin{eqnarray}&\label{AsmpExpnsn}
  (p,q,s)=\sum\limits_{k=0}^{K-1}({p}_k,s_k)(x,t,\xi,\eta)\delta^k+O(\delta^K),\ \delta\to+0,\ K\in \mathbb{N}.
&\end{eqnarray} 
 where the order of approximation, $K$, is as high as one wishes  and the coefficients, ${p}_k,s_k$ are periodic in the fast coordinates, $\xi,\tau$, with periods $\ell_1,\ldots\ell_n,\ell_0$.
 \section{Preliminaries}
\label{ScPrlmnr}
In this section we consider a number of tools for future needs. 
 \subsection{Averaging}
\label{SscAvrgng}
\noindent
Let the skew brackets mean  the averaging over the fast variables, and let operator $M$ do the same -- that is,
\begin{eqnarray}&\label{DfAvrfng}
  (Mf)(x,t)=\langle f \rangle(x,t)\byd\lim\limits_{L\to\infty}\frac{1}{(2L)^{n+1}}\int\limits_{{\mathbb{R}}^{n+1}}f(x,t,\xi,\tau)\,d\xi d\tau
&\end{eqnarray}
In particular, $\langle h \rangle=0\quad \forall x,t$ according the assumption above. Additionally,  we'll  be using the partial (spatial and temporal) averaging operations, that reads
\begin{eqnarray}
 & (M^\xi f)(x,t,\tau)=\langle f \rangle^\xi (x,t,\tau)\byd\lim\limits_{L\to\infty}\frac{1}{(2L)^{n+1}}\int\limits_{{\mathbb{R}}^{n}}f(x,t,\xi,\tau)\,d\xi,&
\label{DfSptAvrgng} \\  
&  (M^\tau f)(x,t,\xi)=\langle f \rangle^\tau (x,t,\xi)\byd\lim\limits_{L\to\infty}\frac{1}{(2L)}\int\limits_{{\mathbb{R}}^{1}}f(x,t,\xi,\tau)\,d\tau.&
\label{DfTmprlAvrgng}
\end{eqnarray}
Evidently,
\begin{eqnarray*}&
M^\xi M^\tau=  M^\tau M^\xi =M \Leftrightarrow\langle\langle f \rangle^\xi \rangle^\tau=\langle\langle f \rangle^\tau \rangle^\xi=\langle f \rangle.
&\end{eqnarray*}
 
Let  $\mathbb{T}^{n+1}=\mathbb{T}^{n+1}_{\xi,\tau}$ be the $n+1$-torus that results from  the factorization of $\mathbb{R}^{n+1}$ over a lattice 
\begin{eqnarray*}
&
\mathbb{Z}_\ell^{n+1}\byd \ell_1\mathbb{Z}\times \ldots\times{\ell_n}\mathbb{Z}\times \ell_{0}\mathbb{Z}.
&
\end{eqnarray*}
Similarly, we identify  n-torus $\mathbb{T}^{n}=\mathbb{T}^{n}_\xi$  with the factor  of  $\mathbb{R}^{n}=\mathbb{R}^{n}_\xi$ over the lattice 
\begin{eqnarray*}
&
\mathbb{Z}_\ell^{n}\byd \ell_1\mathbb{Z}\times \ldots\times{\ell_n}\mathbb{Z}.
&\end{eqnarray*}
Given the  assumptions on the periodicity above,  we identify the fast coordinates, $\xi,\tau$ with  the coordinates on  torus $\mathbb{T}^{n+1}$ and the averaging \eqref{DfAvrfng} with the averaging over this torus. Similarly, we identify the partial mean values \eqref{DfSptAvrgng} and  \eqref{DfTmprlAvrgng} with averaging over tori of $\mathbb{T}_\xi^n$ and $\mathbb{T}^1_\tau$ (the latter is just the circumference).
\subsection{Operators and projectors}
\label{SscDfOpr}
From here on,  using the notation of fast variables as the lower indices at the notations of a dependent variable indicates the partial derivative of this dependent variables in the corresponding fast variables. Similarly, using the notation of fast variables as the lower indices for the notation of a differential operation means that this operation has to be performed by differentiation in the fast variables. 

Let $\efr\in \mathrm{C}^{\infty}\left(\mathbb{T}_{\xi,\tau}^{n+1}\right)$ be an everywhere positive function. Throughout this article,  the notation is due  as follows
\begin{eqnarray}
&
  \rho^{-1}=\langle \efr\rangle^\xi.
&
\label{Dfrho}\\
&\widetilde{\nabla}=\nabla_\xi, \quad\widetilde{\nabla}_{\efr}=\efr \widetilde{\nabla} \frac{1}{\efr},\quad \nabla_{\efr} =\efr\nabla \frac{1}{\efr},\quad \nabla\cdot=\Div,\quad \widetilde{\nabla}\cdot =\Div_\xi,
&
\label{Ntn1}\\
&\mathcal{L}=\widetilde{\nabla}\cdot\widetilde{\nabla}_{\efr},
\label{Ntn2}\\
& \mathcal{H}=\pr_\tau-\Delta_\xi,\quad H=\pr_t-(\widetilde{\nabla}\cdot \nabla +\nabla\cdot\widetilde{\nabla})&
\label{Ntn3}
\end{eqnarray}
Define
\begin{eqnarray}
&
 \mathrm{H}^{s,2}\byd \{u\in{\mathrm{L}}_2\left(\mathbb{T}^{n+1}\right):\, \sum\limits_{} (k^{2s}+m^4)|\hat{u}_{km}|^2<\infty\quad  (k,m)\in \mathbb{Z}_\ell^{n+1}\}
 &
 \label{DfDomH}\\
 &
 \hat{u}_{k,m}=\langle u,\mathrm{e}_{-k,-m}\rangle ,\quad \mathrm{e}_{k,m}={\exp(i(k\tau+m\xi))}{\sqrt{\ell_0\ldots \ell_n}},\quad (k,m)\in \mathbb{Z}_\ell^{n+1}
 &
 \nonumber\\
\end{eqnarray}
Further, define operators $\mathcal{H}:{\mathrm{L}}_2\left(\mathbb{T}^{n+1}\right)\to {\mathrm{L}}_2\left(\mathbb{T}^{n+1}\right)$ and $\mathcal{L}:{\mathrm{L}}_2\left(\mathbb{T}^{n+1}\right)\to {\mathrm{L}}_2\left(\mathbb{T}^{n+1}\right)$ on domains  $\mathrm{H}^{s,2}\subset {\mathrm{L}}_2\left(\mathbb{T}^{n+1}\right)$, $s=1,0$, correspondingly,  by the expressions specified in \eqref{Ntn2} and \eqref{Ntn3}. By this definitions,  operators $\mathcal{H}$ and $\mathcal{L}$ are densely defined and closed.

Operator $\mathcal{H}$ is nothing else than the heat operator on torus, regarding which it is well-known that its kernel
 is the  spectral subspace spanned by the identical eigenfunction.  In other words,  $\ker \mathcal{H}$ consists of the identical functions and there are no  generalized eigenfunctions.  The corresponding spectral projector reads 
\begin{eqnarray*}&
1_*\otimes 1_*,\quad 1_* \sqrt{\ell_0\ldots \ell_n}\stackrel{\mathrm{def}}{\equiv} 1,\quad (a\otimes b)c \byd(c,b)a.
&\end{eqnarray*}
and it is worth noticing that
\begin{eqnarray*}&
1_*\otimes 1_*=M.
&\end{eqnarray*}  
Correspondingly,
\begin{eqnarray*}&
\mathrm{range}\left(I-1_*\otimes 1_*\right)=\{u\in{\mathrm{L}}_2\left(\mathbb{T}^{n+1}\right):\, \langle u \rangle=0 \}\byd \widetilde{\mathrm{L}}_2\left(\mathbb{T}^{n+1}\right).
&\end{eqnarray*}
There exists the right inverse operator, $\mathcal{H}^{-1}$, and it is bounded $\mathcal{H}^{-1}:\widetilde{\mathrm{L}}_2\left(\mathbb{T}^{n+1}\right)\to \widetilde{\mathrm{L}}_2\left(\mathbb{T}^{n+1}\right)$. 
Hence, for every $v\in {\mathrm{L}}_2\left(\mathbb{T}^{n+1}\right)$ such that $\langle v\rangle=0$ the general solution to equation $\mathcal{H}u=v$ on torus $\mathbb{T}^{n+1}$ reads
\begin{eqnarray*}&
u=\langle u\rangle +\mathcal{H}^{-1}v, \quad \langle\mathcal{H}^{-1}v\rangle=0,
&\end{eqnarray*} 
where $\langle u\rangle$ is free parameter. We need the similar information regarding operator $\mathcal{L}$.  For this purpose, define bounded projector
$P:{\mathrm{L}}_2\left(\mathbb{T}^{n+1}\right)\to {\mathrm{L}}_2\left(\mathbb{T}^{n+1}\right)$ by setting
\begin{eqnarray}&\label{DfP-Prjtr}
P:u\to \langle u\rangle^\xi \efr_* \quad  \efr_*=\rho\efr, 
&\end{eqnarray}
\begin{prpsn}\label{PrpsOnPrjctrs} 
\textsf{\emph{Projector $P$ is the $\mathcal{L}$-invariant  projector onto  the kernel of operator $\mathcal{L}:{\mathrm{L}}_2\left(\mathbb{T}^{n+1}\right)\to {\mathrm{L}}_2\left(\mathbb{T}^{n+1}\right)$ -- that is,
\begin{eqnarray*}&
P\mathcal{L}-\mathcal{L}P=0.
&\end{eqnarray*}
}}
\end{prpsn}
{\textbf{Proof.}} For every fix $\sigma\in \mathbb{R}$, let $\breve{\mathcal{L}}=\breve{\mathcal{L}}(\sigma):{\mathrm{L}}_2\left(\mathbb{T}^{n}\right)\to {\mathrm{L}}_2\left(\mathbb{T}^{n}\right)$ be  the operator defined by mapping $u\mapsto \mathcal{L}|_{\tau=\sigma} u$ on  domain 
\begin{eqnarray*}
&
\breve{\mathrm{H}}^2\byd \{u\in{\mathrm{L}}_2\left(\mathbb{T}^{n}\right):\, \sum m^4|\hat{u}_{m}|^2<\infty,\quad m\in \mathbb{Z}_\ell^{n}\}.
&\end{eqnarray*} 
By this definition, operator $\breve{\mathcal{L}}$ becomes strictly elliptic with compact resolvent (if we remind that coefficient $\efr$ is strictly positive by assumption). Direct inspection shows that $\ker\breve{\mathcal{L}}(\tau)$ includes function $\efr(\cdot,\tau)$. 
Let's show that this function  spans $\ker\breve{\mathcal{L}}$. The conjugated operator, $\breve{\mathcal{L}}^*(\tau)$,  acts as follows
\begin{eqnarray*}&
\breve{\mathcal{L}}^*(\tau) u\mapsto \efr^{-1}\widetilde{\nabla}\cdot (\efr \widetilde{\nabla } u)= {\Delta }_\xi u+\widetilde{\nabla } u\cdot \widetilde{\nabla }\ln \efr,\quad \efr=\efr(\cdot,\tau).
&\end{eqnarray*}
Subspace $\ker\breve{\mathcal{L}}^*(\tau)$ consists of solutions to the second order elliptic equation that obeys the strong maximum principle, from which it follows that  every such a solution is an identical function. Hence, $\dim\ker\breve{\mathcal{L}}(\tau)=1$ by the compactness argument.  

Further, define 
\begin{eqnarray*}&
\breve{P}(\tau):u\to \langle u(\cdot,\tau)\rangle^\xi \efr_*(\cdot,\tau).  
&\end{eqnarray*}
A straightforward inspection shows that this is the projector on $\ker\breve{\mathcal{L}}(\tau)$ that commutes with $\breve{\mathcal{L}}^*(\tau)$. Finally, we put $Pu=\breve{P}(\tau)u(\cdot,\tau)$ for $u\in \mathrm{C}(\mathbb{T}^{n+1})$ and then extend it by continuity. 
\\
$\square$
\\
\textbf{Remark 1.} The chain of generalized eigenfunctions linked to $\ker\breve{\mathcal{L}}(\tau)$ (if any) must begin  with the solution to equation $\breve{\mathcal{L}}(\tau)w =\efr(\cdot,\tau)$, but there is a solvability condition to hold. It reads $\langle \efr(\cdot,\tau)\rangle=0$, but   this equality  is never true  by the assumption on positiveness of function $\efr$ above. Thus,   $\ker\breve{\mathcal{L}}(\tau)$ is the spectral subspace in the sense it does not give rise to any generalized eigenfunctions.

It is worth noting that
 \begin{eqnarray}&\label{EqPandMAndMxi}
P(I-M^\xi+M)(\cdot)=(M(\cdot))\efr_*,\ MP=M,\ M^\xi P= M^\xi,\  MQ=M^\xi Q=0,\ 
&\end{eqnarray}
These equalities  follow directly from the definitions.  Furthermore, we have the decomposition of unity by the following triad 
 \begin{eqnarray*}&
 Q=I-P,\quad P_{1}= P(I-M^\xi+M)P,\quad Q_{1}=P(M^\xi-M)P,
 &\end{eqnarray*}
and it induces the decomposition of the ambient space into  the direct sum. This holds true  since every pairwise product is zero in the above triad. The pointwise version of this  decomposition reads
\begin{eqnarray*}&
 u=\langle u\rangle\efr_*+\left(\langle u\rangle^\xi-\langle u\rangle\right)\efr_*+(u-\langle u\rangle^\xi\efr_*)=P_1u+Q_1u+Qu.
 &\end{eqnarray*}
Note in passing, that $M(Q_1+Q)=0$, but $Q_1+Q\neq I-M$. 

Further, define mapping between smooth vector fields in the fast variables on $\mathbb{T}^n=\mathbb{T}^n_\xi$ as follows
\begin{eqnarray}&\label{DfPrjctVctr}
  \mathcal{P}:\vb \mapsto \widetilde{\nabla}_{\efr}\mathcal{L}^{-1}\widetilde{\nabla}\cdot \vb.
&\end{eqnarray}
From here on and until a special notice, by averaging, we mean purely spatial averaging, and we skip  index $\xi$ in its notation.
\begin{prpsn}\label{PrpOnPrjctVctrs}
\textsf{\emph{Mapping \eqref{DfPrjctVctr}  allows for extending to  the bounded projector ${\mathrm{L}}_2\left(\mathbb{T}^{n}\right)\to {\mathrm{L}}_2\left(\mathbb{T}^{n}\right)$ such that
\begin{eqnarray*}&\mathrm{range\,}{\mathcal{P}}=\{\vb=\widetilde{\nabla}_{\efr}\phi,\ \phi\in \mathrm{H}^{1}\left(\mathbb{T}^{n}\right)\},\quad
\ker\mathcal{P}=\{\vb:\widetilde{\nabla}\cdot\vb=0\}.&\end{eqnarray*}
}}
\end{prpsn}
{\bf Proof.} The boundedness follow from the common elliptic estimates for equations in the divergent form. Further, 
\begin{eqnarray*}&
\mathcal{P}^2=\widetilde{\nabla}_{\efr}\mathcal{L}^{-1}\widetilde{\nabla}\cdot\widetilde{\nabla}_{\efr}\mathcal{L}^{-1}\widetilde{\nabla}\cdot=
\widetilde{\nabla}_{\efr}\mathcal{L}^{-1}\mathcal{L}\mathcal{L}^{-1}\widetilde{\nabla}\cdot=
\widetilde{\nabla}_{\efr}\mathcal{L}^{-1}\widetilde{\nabla}\cdot=\mathcal{P}.
&\end{eqnarray*}
Finally, $\mathcal{L}^{-1}\mathcal{L}\phi=\phi-\langle\phi\rangle\ve_*$, and 
\begin{eqnarray*}&
\mathcal{P}\widetilde{\nabla}_{\efr}\phi
=\widetilde{\nabla}_{\efr}\mathcal{L}^{-1}\widetilde{\nabla}\cdot\widetilde{\nabla}_{\efr}\phi=\widetilde{\nabla}_{\efr}\mathcal{L}^{-1}\mathcal{L}\phi=
\widetilde{\nabla}_{\efr}  (\phi-\langle\phi\rangle\ve_*)=\widetilde{\nabla}_{\efr}\phi.\quad \square
&\end{eqnarray*}
{\bf Remark 2.} Projector $\mathcal{P}$ is not orthogonal relative to the standard Euclidean structure in space ${\mathrm{L}}^v_2\left(\mathbb{T}^{n}\right)$ (superscript means that this notation is for the space of vector fields). However, it gets orthogonal in metric 
\begin{eqnarray}&\label{DfE-mtrc}
\|\vf\|_\efr^2\byd \langle |\vf|^2 /\efr\rangle
&\end{eqnarray}
Averaging $\mathcal{P}$-projections of the special vector fields that read $\efr \va$, $\va=\const$, results in a linear operator   $\mathcal{M}:\mathbb{R}^n\to \mathbb{R}^n$ that acts as follows
\begin{eqnarray}&\label{DfMtrxM}
\forall\, b\in \mathbb{R}^n\  b\cdot \mathcal{M}a=(\efr\vb,\mathcal{P}(\efr \va))_{\efr}\quad \va=J a,\ \vb=Jb,
&\end{eqnarray}
where $a \cdot b $ -- standard Euclidean structure on $\mathbb{R}^n$, $J$ is the embedding $\mathbb{R}^n\to \mathrm{L}^v_{2,e}\left(\mathbb{T}^{n}\right)$ induced by Cartesian coordinates, and  $\mathrm{L}^v_{2,e}\left(\mathbb{T}^{n}\right)$ is the vector fields space with Euclidean metric \eqref{DfE-mtrc}.
\begin{prpsn}\label{PrpOnM-mtrx}
\textsf{\emph{Operator  $\mathcal{M}$ is symmetric, its quadratic form is non-negative, and its eigenvalues are strictly smaller than $\langle \efr \rangle$.}}
\end{prpsn}
\textbf{Proof.} Since the projector $\mathcal{P}$ is orthogonal in metric  \eqref{DfE-mtrc}, it is symmetric and its norm is less than unity. Hence, bilinear form \eqref{DfMtrxM} is symmetric. Besides, $(\efr\va,\mathcal{P}(\efr \va))_{\efr}\ge 0$, so that  
$
a\cdot\mathcal{M}a \ge 0
$
too. Finally,
\begin{eqnarray*}&|b\cdot \mathcal{M}a|\le \|\efr\va\|_\efr\|\efr\vb\|_\efr=\langle\efr\rangle|a||b|,&\end{eqnarray*}
and this yields the upper bound for eigenvalues. Assume this bound is attained for some $a$. Then $\efr\va\in \mathrm{range}\, \mathcal{P}$ -- that is, $\efr\va=\widetilde{\nabla} _\efr\phi $, but this is impossible since a constant vector field on torus has no single-valued potential.  
 $\square$
 \\
\textbf{ Remark 3. } When function $\efr$ possesses some  translational invariance,  form $a\cdot\mathcal{M}a$  gets degenerated. Its null-subspace  in $\mathbb{R}^n$ is  the one identified with the subspace of those translations  that leave function $\efr$ invariant. (Here we mean the translations in the fast spatial variables, $\xi_1,\ldots,\xi_n$, while all the others are frozen.)
 \\
\textbf{ Remark 4. } For every  fixed  $\va$, bilinear form \eqref{DfMtrxM} determines a directional differentiation, $v$,  in slow variables by expression 
\begin{eqnarray*}&
d\eta(v)=(\efr\nabla\eta,\mathcal{P}(\efr \va))_{\efr}, 
&\end{eqnarray*} 
or, in some coordinates,
$v_j=(\efr\nabla x_j,\mathcal{P}(\efr \va))_{\efr}$.
Hence, relative to  Cartesian coordinates with basis $b_1,b_2,..b_n$, 
\begin{eqnarray*}&
v_j=(\efr\vb_j ,\mathcal{P}(\efr \va))_{\efr}, j=1..n,
&\end{eqnarray*}
Similarly, the entries of matrix of operator $\mathcal{M}$ read
\begin{eqnarray*}&
\mu_{ij}=b_i\cdot \mathcal{M}b_j=(\efr b_i,\mathcal{P}(\efr b_j))_{\efr},\ i,j=1..n
&\end{eqnarray*}
At this point and from it on, we restore the notations, which we had been  introducing earlier  for  distinguishing between the spatial, temporal and total averages.
\section{ Asymptotic}
\label{ScAsmpt}
\noindent
Here we get back to the system consisting of equations \eqref{EqPrdTrns}, \eqref{EqFlx} and \eqref{EqPryDnst}, where the external signal, $h$, is a general short wave  given by expression \eqref{ShrtWvSgnl}.   

The system under consideration written with the use of the fast variables reads as
\begin{eqnarray}
&p_\tau+\widetilde{\nabla}\cdot \vq =\delta\left(pf(p,s)-p_{t}-{\nabla}\cdot\vq\right),&
\label{EqPrdTrnsFstVrbl0}\\
&p\widetilde{\nabla}\left(\chi s+\kappa h -\mu\ln p\right)= \delta\left(\vq-p\nabla\left(\chi s+\kappa h -\mu\ln p\right)\right),&
\label{EqFlxFstVrbl0}\\
&(\mathcal{H}+\delta H){s}=\delta\left(sg(p,s)+\delta\Delta s\right),&
\label{EqPryDnstFstVrbl0}
 \end{eqnarray}
 The following ansatz 
\begin{eqnarray}&\label{ansatz}
p(x,t,\xi,\tau)=r(x,t,\xi,\tau)\efr(x,t,\xi,\tau),\quad \efr\byd\exp(\kappa h/\mu),\quad h=h(x,t,\xi,\tau) 
&\end{eqnarray} 
brings the system into the form
\begin{eqnarray}
&p_\tau+\widetilde{\nabla}\cdot \vq =\delta\left(pf(p,s)-p_{t}-{\nabla}\cdot\vq\right),\quad p=r\efr,&
\label{EqPrdTrnsFstVrbl}\\
&r\widetilde{\nabla}\left(\chi s-\mu\ln r\right)= \delta\left(\vq/\efr-r\nabla\left(\chi s -\mu\ln r\right)\right),&
\label{EqFlxFstVrbl}\\
&(\mathcal{H}+\delta H){s}=\delta\left(sg(p,s)+\delta\Delta s\right),\quad p=r\efr.&
\label{EqPryDnstFstVrbl}
 \end{eqnarray} 
Replacing the unknowns in equations \eqref{EqPrdTrnsFstVrbl}-\eqref{EqPryDnstFstVrbl} by expansions \eqref{AsmpExpnsn} with $p_k=r_k\efr $ 
leads to a chain of equations  for the   coefficients $r_k,q_k,s_k$, $k=0,1,2,\ldots$, and it reads as 
 \begin{eqnarray}
&p_{k\tau}+\widetilde{\nabla}\cdot \vq_k =-p_{k-1,t}-{\nabla}\cdot\vq_{k-1}+\sum\limits_{j=0..k-1}p_jf_{k-1-j},\quad p_k=r_k\efr.&,
\label{EqPrdTrnsCfcnt}\\
&-\mu \widetilde{\nabla} r_k+\chi \sum\limits_{j=0..k}r_j\widetilde{\nabla} s_{k-j}=\vq_{k-1}/\efr+\mu \nabla r_{k-1}-\chi \sum\limits_{j=0..k-1}r_j\nabla s_{k-1-j},&
\label{EqFlxCfcnt}\\
&\mathcal{H}{s}_k+H s_{k-1} = \Delta s_{k-2}+\sum\limits_{j=0..k-1}s_jg_{k-1-j},\quad p_k=r_k\efr.&
\label{EqPryDnstCfcnt}
 \end{eqnarray}
Here and  in the sequel  all the quantities with negative indices are equal to zero in every equation, and  coefficients $f_k,g_k$ result from expanding the   reaction terms as follows  
\begin{eqnarray}&\label{DfFm-Gm}
  f(p,s)=\sum\limits_{m=0}f_m\delta^m,\  f_0=f(p_0,s_0),\quad  g(p,s)=\sum\limits_{m=0}g_m\delta^m,\  g_0=g(p_0,s_0),
&\end{eqnarray}
where we have omitted  the explicit expressions  for the exposition compactness.

We recall that equations \eqref{EqPrdTrnsCfcnt}-\eqref{EqPryDnstCfcnt} have to be solved in the functions, which are well-defined on the torus of the fast variables, $\mathbb{T}_{\xi,\tau}^{n+1}$. 

We remind three projectorsintroduced in  Subsection~\ref{SscDfOpr}. These are 
\beq
P_1=P(I-M^\xi+M)P,\quad Q_1=P(M^\xi-M)P,\quad Q=I-P
\eeq
The  notations for the correspondent projections of coefficients $p_k,\ k=0,1,2,\ldots,$ read
\begin{eqnarray*}&
\pht^\bullet_k,\quad\pht^\circ_k,\quad \pch_k.
&\end{eqnarray*}
Clearly,
\begin{eqnarray*}&
\pht^\bullet_k+\pht^\circ_k=\pht_k.
&\end{eqnarray*}
Additionally, we set 
\begin{eqnarray*}&
\pch^\circ_k\byd \pht^\circ_k+\pch_k.
&\end{eqnarray*}
Also, for the averages and for the deviations from them,  the notations is as follows
\begin{eqnarray*}
&
\widetilde{p}_k=(I-M)p_k,\quad\widetilde{p}^\tau_k=(M^\xi-M) p_k,\quad \widetilde{s}_k=(I-M)s_k,
&
\\
&
\bar{p}_k=Mp_k,\quad  \bar{s}_k=Ms_k,
\quad \bpx_k=M^\xi_k.
&
\end{eqnarray*}
Note that 
\begin{eqnarray*}&
M\pht^\bullet_k=M\pht_k=Mp_k,, \quad M^\xi\pht_k=M^\xi p_k,\quad  M^\xi\pch_k=0,\quad M\pch^\circ_k=0,
&\end{eqnarray*}
but  $M^\xi\pch^\circ_k\neq0$, $\pch^\circ_k\neq (I-M)p_k$ generally. 
\subsection{Leading approximation}\label{SScLdngApprx}
\noindent
For $k=0$, equations \eqref{EqPrdTrnsCfcnt}-\eqref{EqPryDnstCfcnt} read
\begin{eqnarray}
&p_{0\tau}+\widetilde{\nabla}\cdot\vq_0 =0,\quad p_0=r_0\mathfrak{e}&,
\label{EqPrdCfcOrd0}\\
&\widetilde{\nabla}\left(\chi s_0 -\mu\ln r_0\right)=0,&
\label{EqFlxCfcOrd0}\\
&\mathcal{H}{s}_0=0.&
\label{EqPryCfcOrd0}
 \end{eqnarray}
The solution to equations \eqref{EqPryCfcOrd0} reads
\begin{eqnarray}&\label{S0}
s_0=\bar {s}(x,t).
&\end{eqnarray}
Hence,
\begin{eqnarray*}&
r_0=r_0(x,t,\tau),
&\end{eqnarray*}
by equation \eqref{EqFlxCfcOrd0}, and equation \eqref{EqPrdCfcOrd0} entails the following
\begin{eqnarray*}&
\left(r_0\langle \efr \rangle^\xi\right)_\tau=0.
&\end{eqnarray*}
By the last equation,  
\begin{eqnarray}&\label{P0}
p_0=\pht^\bullet_0=\bar{p}(x,t)\efr_*(x,t,\xi,\tau),\ \efr_*=\rho{\efr}
 &\end{eqnarray}

Thus, in the leading approximation, the predators respond to the short wave   external signal by a short wave pattern, $p_0=\bar{p}\efr_*$, that undergoes a slow modulation due to the slowly varying amplitude, $\bar{p}=\bar{p}(x,t)$, which, in turn, controls over  the mean predators density, $Mp_0=\bar{p}$.  The system  that governs the slow modulation of the leading approximation  turns out to be a byproduct of calculating the next approximation, for which we have equations as follows
\begin{eqnarray}
&p_{1\tau}+ \widetilde{\nabla}\cdot\vq_1 +p_{0,t}+\nabla\cdot\vq_{0}= p_0f_{0}.&
\label{EqPrdCfcOrdr=1}\\
&\chi p_0\widetilde{\nabla} s_{1}-\mu \widetilde{\nabla}_{\efr} {p}_1 +\chi p_0\nabla s_{0}-\mu\nabla_{\efr} p_{0}=\vq_{0},&
\label{EqFlxCfcOrdr=1}\\
&\mathcal{H}s_{1}+ s_{0t}=s_0g_{0}.&
\label{EqPryCfcOrdr=1}
 \end{eqnarray}
We have  dropped the zero term,  $p_{1}\widetilde{\nabla} s_{0}$. It's correct  since $s_0=\bar{s}(x,t)$.
\begin{dfntn}\label{DfLdngSlwSst}
\textsf{\emph{By the leading slow system we mean the system of equations relative to the mean values of the leading approximation, $\bar{p}=\bar{p}(x,t)$, $\bar{s}=\bar{s}(x,t)$, that reads
\begin{eqnarray}
&
\bar{s}_{t}=\bar{s}\bar{g}(\bar{p},\bar{s}),\quad \bar{g}(\bar{p},\bar{s})=\langle g(\efr_*\bar{p},\bar{s})\rangle.
&
\label{EqBrS}\\
&
\bar{p}_{t}+\nabla\cdot\left(\bar {p}(\bar{\vc}+\chi \bar{\mathcal{D}}\nabla \bar {s}) - \mu\bar{\mathcal{D}}\nabla \bar{p}\right)=\bar{p}\bar{f}(\bar{p},\bar{s}),\quad \text{where}\ \bar{f}(\bar{p},\bar{s})=\langle \efr_* f(\efr_*\bar{p},\bar{s})\rangle,
&
\label{EqBrP}\\
&
\bar{\vc}\byd \mu\left\langle\left(E-\rho\mathcal{M}\right)\nabla\ln \langle\efr\rangle^\xi \right\rangle^\tau -\langle \widetilde{\nabla}_{\efr}\mathcal{L}^{-1}\efr_{*\tau}\rangle,\quad \bar{\mathcal{D}}\byd E-\left\langle\rho\mathcal{M}\right\rangle^\tau,
&
\label{DfDrftAndVscTnsrLdng}
\end{eqnarray}
and the notation of $\mathcal{M}$ is for the matrix defined  by equality \eqref{DfMtrxM}. 
}}
\end{dfntn} 
\begin{prpsn}\label{PrpOnMtrxD}
\textsf{\emph{ 
Matrix $\bar{\mathcal{D}}$ is symmetric and positively defined. 
}}
\end{prpsn}
\textbf{{Proof}} follows directly from Proposition~\ref{PrpOnM-mtrx}.
\\$\square$\\
\textbf{Remark 5.} By the  Proposition~\ref{PrpOnMtrxD},  differential operator $\nabla\cdot (\bar{\mathcal{D}}\nabla)$ is elliptic. 
\\
\textbf{Remark 6.} From Remark~3, it follows that matrix $\bar{\mathcal{M}}$ can be degenerated. In particular, it is equal to zero if the external signal is off -- that is, if $h\equiv 0$. Then $\bar{\mathcal{D}}=E$.
  \begin{theorem}\label{PrpsLdngSlwdSst}
\textsf{\emph{Let  the mean values of the leading approximation, $\bar{p}=\bar{p}(x,t)$, $\bar{s}=\bar{s}(x,t)$ deliver a solution to the leading slow system. Then the equations of the next approximation, \eqref{EqPrdCfcOrdr=1}-\eqref{EqPryCfcOrdr=1} are compatible and determine the first order corrections to the densities, $p_{1}$  and $s_{1}$ up to their mean values, $\bar{p}_1$ and $\bar{s}_1$. More precisely, the components $\pht^\bullet_1$ and $\bar{s}_1$ remain undetermined, while the other  ones, $\widetilde{s}_1$,  $\pht^\circ_1$ and $\pch_1$ are defined uniquely.
}}
\end{theorem}
{\bf Proof.} Given expressions \eqref{S0} and \eqref{P0},  we transform the  equations  \eqref{EqPrdCfcOrdr=1}-\eqref{EqPryCfcOrdr=1}  as follows
 \begin{eqnarray} 
&p_{1\tau}+ \widetilde{\nabla}\cdot\vq_1 +(\bar{p}\efr_{*})_{t}+\nabla\cdot\vq_{0}=\bar{p}\efr_{*}f(\bar{p}\efr_{*},\bar{s}).&
\label{EqPrdCfcOrdr=1-1}\\
&\chi \bar{p}\efr_*\widetilde{\nabla}\widetilde{s}_{1}-\mu\widetilde{\nabla}_\efr {p}_1 +
\efr_*\left(\bar{p}\nabla(\chi\bar{s}+\mu \ln \langle\efr\rangle^\xi)
-\mu\nabla\bar{p}\right)=\vq_{0},&
\label{EqFlxCfcOrdr=1-1}\\
&\mathcal{H}{s}_1=\bar{s}g(\bar{p}\efr_{*},\bar{s})-\bar{s}_{t}.&
\label{EqPryCfcOrdr=1-1}
 \end{eqnarray}
 Equation   \eqref{EqPryCfcOrdr=1-1} is the  heat equation on the torus, $\mathbb{T}^{n+1}$. Its  compatibility condition  reads
 \begin{eqnarray*}&
 M\left(\bar{s}g(\bar{p}\efr_{*},\bar{s})-\bar{s}_{t}\right)=0,
 &\end{eqnarray*}
 and this is  equation \eqref{EqBrS}. Since this equation holds true by the assumption,  we  conclude that
 \begin{eqnarray*}&
 s_{1}=\bar{s}_1+\widetilde{s}_1,\quad \widetilde{s}_1=\mathcal{H}^{-1}\left(\bar{s}g(\bar{p}\efr_{*},\bar{s})-\bar{s}_{t}\right),
 &\end{eqnarray*}
 where the notation of $\mathcal{H}^{-1}$ is for the right inverse operator, which sends $\ker M$ to itself. Hence, $\bar{s}_1=Ms_{1}$.
 
 Equation  \eqref{EqPrdCfcOrdr=1-1}  also needs a  compatibility condition    to hold, which is  as follows
 \begin{eqnarray}&\label{EqCmptPrdOrd=1}
  \bpx_{1,\tau}=\left\langle (\bar{p}\efr_{*})_{t}+\nabla\cdot\vq_{0}-\bar{p}\efr_{*}f(\bar{p}\efr_{*},\bar{s})\right\rangle^\xi.
 &\end{eqnarray}
 We consider this one as an equation relative to the spatial average of the first correction,  $\bpx_{1}$. Further,  we eliminate vector $\vq_0$ by setting it as follows
 \begin{eqnarray}&\label{EqQ0Dcmps}
\vq_0=-\bar{p}\widetilde{\nabla}_{\efr}\mathcal{L}^{-1}\efr_{*\tau}+(I-\mathcal{P})\vu,\quad \vu\byd\efr_*\left(\bar{p}\nabla(\chi\bar{s}+\mu \ln \langle\efr\rangle^\xi)-\mu\nabla\bar{p}\right)
&\end{eqnarray}
where  we  employ the orthogonal decomposition delivered  by proposition~\ref{PrpOnPrjctVctrs} and remarks aftermath.
 This setting  is consistent with equation \eqref{EqPrdCfcOrd0}. Indeed, $\widetilde{\nabla}\cdot (I-\mathcal{P})=0$, hence,
\begin{eqnarray*}&
\widetilde{\nabla}\cdot\vq_0=-\bar{p}\widetilde{\nabla}\cdot\widetilde{\nabla}_{\efr}\mathcal{L}^{-1}\efr_{*\tau}=-\bar{p}\efr_{*\tau}=-p_{0\tau}.
&\end{eqnarray*}
 Now we aim at finding function ${p}_1$ that  makes setting  \eqref{EqQ0Dcmps} consistent with equation \eqref{EqFlxCfcOrdr=1-1}. This leads to the following equation
 \begin{eqnarray*}&
 \chi \bar{p}\efr_*\widetilde{\nabla}\widetilde{s}_{1}-\mu \widetilde{\nabla}_\efr {p}_1 =-\bar{p}\widetilde{\nabla}_{\efr}\mathcal{L}^{-1}\efr_{*\tau}-\mathcal{P}\vu.
 &\end{eqnarray*}
 Further, we put this as follows
 \begin{eqnarray}&\label{EqCmptFlxOrdr=1}
 \efr\widetilde{\nabla}\left(\chi \rho\bar{p}\widetilde{s}_{1}-\mu {p}_1/\efr +\efr^{-1}\mathcal{L}^{-1}(\widetilde{\nabla}\cdot\vu+\bar{p}\efr_{*\tau}) \right)=0.
&\end{eqnarray}
Here, the common multiplier, $\efr$, is strictly positive by assumption, hence,  we  remove it, and then isolate ${p}_1$ up to projection $\pht_1$.  The result reads
 \begin{eqnarray}
\mu \pch_1=\chi Q {p}_0\widetilde{s}_{1}+\mathcal{L}^{-1}(\widetilde{\nabla}\cdot\vu+\bar{p}\efr_{*\tau})+\mu \pht_1.
\label{EqPrdCfcNoFlxOrd=1}
\end{eqnarray}
Here, the second summand on the right hand side belongs to $\mathrm{range}\,Q$ as  $M^\xi\mathcal{L}^{-1}=0$ and $P\mathcal{L}^{-1}=0$, therefore. 

Thus, we have defined uniquely the projection on  subspace $\mathrm{range}\, Q$, while the projection on subspace  $\mathrm{range}\, P$, $\pht_{1}$,  remains indefinite yet. Let's employ it to make equation \eqref{EqPrdCfcOrdr=1-1} compatible. Its compatibility condition reads 
\beq
 \langle \pht_{1} \rangle^\xi_{\tau}+ \langle (\bar{p}\efr_{*})_{t}+\nabla\cdot\vq_{0}-\bar{p}\efr_{*}f(\bar{p}\efr_{*},\bar{s})\rangle^\xi=0.
\eeq
Isolating $\pht_{1}$ from here also requires the compatibility condition that reads
\beq 
\langle (\bar{p}\efr_{*})_{t}+\nabla\cdot\vq_{0}-\bar{p}\efr_{*}f(\bar{p}\efr_{*},\bar{s})\rangle=0.
\eeq
We eliminate vector $\vq_0$ from  here  by formula \eqref{EqQ0Dcmps} and arrive at the equation as follows
\begin{eqnarray*}
 &
\bar{p}_t+\nabla\cdot\langle\bar{p}\widetilde{\nabla}_\efr\mathcal{L}^{-1}\efr_{*\tau}+\rho(I-\mathcal{P})\vu_1\rangle=
\langle\bar{p}\efr_{*}f(\bar{p}\efr_{*},\bar{s})\rangle,
 &
 \\
 &
 \vu_1\byd\vu /\rho=\chi \bar{p}\efr\nabla \bar{s}-\mu\efr\nabla\bar{p}
+\mu \bar{p}\efr \nabla\ln \langle\efr\rangle.
 &
\end{eqnarray*}
From definition \eqref{DfMtrxM}, it follows that this  equation is equivalent to equation \eqref{EqBrP}, which, in turn, holds  by assumption. Hence, we can determine component $\pht_{1}$ up to a solution of the following  homogeneous equation 
\beq 
\pr_\tau M^\xi u=0.
\eeq
 Since every such a solution belongs to $\mathrm{range}\, P_1$, we conclude that we are able to determine the projection on $\mathrm{range}\,Q_1$, $\pht^\circ_1$, uniquely. The corresponding formula reads
\beq
\pht^\circ_1=P\pr_\tau^{-1} \langle (\bar{p}\efr_{*})_{t}+\nabla\cdot\vq_{0}-\bar{p}\efr_{*}f(\bar{p}\efr_{*},\bar{s})\rangle^\xi,
\eeq
where the notation of $\pr^{-1}$ is for the right inverse operator that sends $\mathrm{range}\,(I-M^\tau)$ to itself, while vector $\vq_0$ is specified by formula \eqref{EqQ0Dcmps}. Thus, we have determined the first order corrections, $p_1$ and $s_1$, up to projections $\pht^\bullet_1$ and $\bar{s}$, and this completes the proof. 
\\$\square$\\
{\bf Remark 7.}  Note that formulation of Theorem~\ref{PrpsLdngSlwdSst} presumes that  the uniquely defined components are  independent of the ones that remain undetermined in the claimed decomposition.  At the same time, knowing the  averages of the next corrections removes this uncertainty completely, and vice versa.
\subsection{Examples}\label{Ssc1Dmn}
\noindent
In one spatial dimension,  the diffusivity tensor, $\overline{\mathcal{D}}$  and  the drift, $\overline{\vc}$,  get   scalar in effect. We calculate them by formulae \eqref{DfDrftAndVscTnsrLdng}, where we put   
\begin{eqnarray*}
    &
\widetilde{\nabla}\cdot=\widetilde{\nabla}=\pr_\xi,\quad \widetilde{\nabla}_{\efr}=\efr\pr_\xi(\cdot/\efr),\quad  \mathcal{L}(\cdot)=\pr_\xi(\efr\pr_\xi(\cdot/\efr)).
    &
\end{eqnarray*}
Hence,
\begin{eqnarray*}&
\mathcal{M}=\langle\efr\pr_\xi(\efr^{-1}\mathcal{L}^{-1}\pr_\xi\efr)) \rangle^\xi=-\langle\efr^{-1}\pr_\xi\efr\mathcal{L}^{-1}\pr_\xi \efr\rangle^\xi=-\langle\efr^{-1} u \mathcal{L} u\rangle^\xi=\langle\efr(\pr_\xi (u/\efr))^2\rangle^\xi.
&\end{eqnarray*}
Calculating $\mathcal{L}^{-1}\pr_\xi \efr$ reduces to solving ODE
\begin{eqnarray*}&
\pr_\xi(\efr \pr_\xi(u/\efr))=\pr_\xi \efr.
&\end{eqnarray*}
Integrating gives 
\begin{eqnarray*}&
\pr_\xi(u/\efr)=1 -\frac{1}{\efr\langle\efr^{-1}\rangle^\xi}.
&\end{eqnarray*}
Hence,
\begin{eqnarray*}&
\langle\efr(\pr_\xi (u/\efr))^2\rangle^\xi=\left\langle\efr\left(1 -\frac{1}{\efr\langle\efr^{-1}\rangle^\xi}\right)^2\right\rangle^\xi=\langle\efr\rangle^\xi-\frac{1}{\langle\efr^{-1}\rangle^\xi},
&\end{eqnarray*}
and we arrive at the following equality 
\begin{eqnarray*}&
\overline{\mathcal{D}}=1-\langle\rho\mathcal{M}\rangle^\tau=\left\langle \frac{1}{\langle\efr\rangle^\xi\langle\efr^{-1}\rangle^\xi}\right\rangle^\tau.
&\end{eqnarray*}
Further, we pass to the drift.  Let $\bar{c}$ denote the scalar identified with the drift velocity, $\bar{\vc}$. We address its component  that reads 
\begin{eqnarray*}&
-\langle\efr\pr_\xi(\efr^{-1}\mathcal{L}^{-1} \efr_{*\tau})) \rangle.
&\end{eqnarray*} 
Let $w=\pr^{-1}_\xi \efr_{*\tau}$, where we employ the right inverse operator, such that $M^{\xi}\pr_\xi^{-1}=0\ \forall\,\tau$. Then calculating $v=\mathcal{L}^{-1} \efr_{*\tau}$ reduces to solving ODE 
\begin{eqnarray*}&
\pr_\xi(\efr \pr_\xi(v/\efr))=\pr_\xi w.
&\end{eqnarray*}
Hence, 
\begin{eqnarray*}&
\pr_\xi(v/\efr)=(w+C)/\efr,
&\end{eqnarray*}
and
\begin{eqnarray*}&
C=-\frac{1}{\langle\efr^{-1}\rangle^\xi}\left\langle\frac{w}{\efr}\right\rangle^\xi.
&\end{eqnarray*}
So, the drift component under consideration  simplifies  to
\begin{eqnarray*}&
\bar{c}_1=-\langle\efr\pr_\xi(v/\efr) \rangle=-C
&\end{eqnarray*}
Thus, in one spatial dimension, the leading slow system  includes equation \eqref{EqBrS} as is, while equation \eqref{EqBrP} takes simpler form as follows  
\begin{eqnarray}
&
\bar{p}_{t}+\left(\bar {p}(\bar{c}+\bar{\chi} \nabla \bar {s}) - \bar{\mu}\bar{p}_x)\right)_x=\bar{p}\bar{f}(\bar{p},\bar{s}),\ \ \quad \text{where}\quad \bar{f}(\bar{p},\bar{s})=\langle \efr_* f(\efr_*\bar{p},\bar{s})\rangle,
&
\label{EqBrP1D}\\
&
\bar{\chi}=\left\langle\frac{\chi }{\langle\efr\rangle^\xi\langle\efr^{-1}\rangle^\xi}\right\rangle^\tau, \quad \bar{\mu}=\left\langle\frac{\mu}{\langle\efr\rangle^\xi\langle\efr^{-1}\rangle^\xi}\right\rangle^\tau,\bar{c}=\left\langle \bar{\mu}\pr_x\ln \langle\efr\rangle^\xi +\frac{1}{\langle\efr^{-1}\rangle^\xi}\left\langle\frac{w}{\efr}\right\rangle^\xi\right\rangle^{\tau},
&
\label{DfDrftAndVscTnsrLdng1D}
\end{eqnarray}
and $\pr_\xi w=\efr_{*\tau}$, $\langle w\rangle^\xi=0\ \forall\ \tau$.

More simplifications become feasible if we  assume the signal to have the form  of a short traveling wave  -- that is, if we set
\beq
h=h(\eta),\ \eta=\xi-c\tau.
\eeq
Then spatial and temporal averages are equal, and $\pr_\tau=-c\pr\eta=\pr_\xi$. In particular, 
$
w=-c(\efr_*-1).
$
In this way, formulae \eqref{DfDrftAndVscTnsrLdng1D} gets simpler and now read
\begin{eqnarray*}&
\bar{\chi}=\frac{\chi }{\langle\efr\rangle^\xi \langle\efr^{-1}\rangle^\xi }, \quad \bar{\mu}=\frac{\mu}{\langle\efr\rangle^\xi \langle\efr^{-1}\rangle^\xi },
\bar{c}=c\left(1-\frac{1}{\langle\efr\rangle^\xi \langle\efr^{-1}\rangle^\xi} \right)+ \bar{\mu}\pr_x\ln\langle\efr\rangle^\xi.
&\end{eqnarray*}
Thus, we have arrived to a particular form of the leading slow system,  which is exactly the one that  Morgulis and Malal \cite{AM3} had already derived in the case of one spatial dimension assuming additionally that the external signal and the corresponding  solution represent the traveling waves spreading over the fast variables.  

Further, in multiple dimensions, suppose the signal  to read as follows
\beq
h=h(\eta),\ \eta=\theta\cdot\xi-c\tau,\quad  |\theta|=1.
\eeq
Then, simply repeating the calculations above, we arrive at a multidimensional version of the leading slow system, where 
\begin{eqnarray*}&
\overline{\mathcal{D}}=1-\theta\otimes\theta+ \frac{\theta\otimes\theta}{\langle\efr\rangle^\xi\langle\efr^{-1}\rangle^\xi}, \quad \bar{\vc}=\bar{c}\vtheta,
&\end{eqnarray*}
and we have got specified parameter  $\bar{c}$ above. Of course, we can consider a slowly modulated traveling wave,  $h=h(x,t,\eta)$. It  can be addressed in the same way.
 
 It is worth noticing, that 
\beq 
\langle\efr\rangle^\xi\langle\efr^{-1}\rangle^\xi>1.
\eeq
Indeed, $\langle\efr\rangle^\xi\ge \exp(\varkappa\langle h\rangle^\xi/\mu)$ and  $\langle\efr^{-1}\rangle^\xi\ge \exp(-\varkappa\langle h\rangle^\xi/\mu)$ by convexity. 
To make the formulae above fully explicit, suppose, for example,  $h(x,t,\xi,\tau)=a(x,t,\tau)\cos(\xi+\xi_0)$, where $\xi_0$ is a constant phase. Then 
\begin{eqnarray*}
&
\langle\efr\rangle^\xi= \langle\efr^{-1}\rangle^\xi=I_0(a).
&
\end{eqnarray*}
where $I_0$ is the second kind Bessel function of order 0. Thus, perceiving an external signal on average can lead to pressing down both the diffusivity and the sensitivity to the other signals. 
\subsection{Higher  approximations}
\label{SscHghOrdrApprx}
\noindent
The system governing the corrections of order $k$, $k=2,\ldots$, reads
  \begin{eqnarray}
&p_{k\tau}+ \widetilde{\nabla}\cdot\vq_k +p_{k-1,t}+\nabla\cdot\vq_{k-1}=\sum\limits_{j=0..k-1}p_jf_{k-1-j}.&
\label{EqPrdCfcOrdr>0}\\
&-\mu (\widetilde{\nabla}_{\efr} {p}_k + \nabla_{\efr} p_{k-1}) +\chi\sum\limits_{j=0..k-1}p_j(\widetilde{\nabla} {s}_{k-j}+\nabla s_{k-1-j})=\vq_{k-1},&
\label{EqFlxCfcOrdr>0}\\
&\mathcal{H}{s}_k+H s_{k-1}-\Delta s_{k-2}=\sum\limits_{j=0..k-1}s_jg_{k-1-j}.&
\label{EqPryCfcOrdr>0}
 \end{eqnarray}
Since $s_0$ is purely slow function, we have dropped the term involving $\widetilde{\nabla} s_0$ from equation \eqref{EqFlxCfcOrdr>0}. 
It's useful to see that actually
\begin{eqnarray}
&\vq_{k-1}=-\mu (\widetilde{\nabla}_{\efr} \pch_k + \nabla_{\efr} p_{k-1}) +\chi\sum\limits_{j=0..k-1}p_j(\widetilde{\nabla} \widetilde {s}_{k-j}+\nabla s_{k-1-j})&
\label{EqFlxCfcOrdr>0-1}
 \end{eqnarray}
 provided that equation \eqref{EqFlxCfcOrdr>0} holds for some functions $p_k$, $s_{k-j}$.
 
 For  compactness of layout, we'll be  using  a notation of $\mathrm{Op}_m$, $m=1,2,3\dots$, in place of  an  operation (maybe, differential) that involve only those coefficients  $p_k,\vq_k,s_k$ that are indexed by $k\le m$. 
Further,  we put
\begin{eqnarray}
&
a_{i}({\zeta}_1,{\zeta}_2)\byd a_{ij}{\zeta}_j,\ i,j=1,2,\quad  a_{ij}\byd F_{iz_j}\left|_{z_1=p_0,z_2={s}_0}\right.,\quad\text{where}
&
\label{EqDf-a-LnFn}\\
& 
\ F_1=z_1f(z_1,z_2),\ F_2=z_2g(z_1,z_2),
&
\nonumber
\end{eqnarray} 
Then
\begin{eqnarray}
&
\sum\limits_{j=0..m}p_jf_{m-j}-a_1(p_{m},s_{m})=\sum\limits_{j=1..m-1}p_jf_{m-j}=Op_{m-1},
&
\nonumber\\
&
\sum\limits_{j=0..m}s_jg_{m-j}-a_{2}(p_{m},s_{m})=\sum\limits_{j=1..m-1}s_jg_{m-j}=Op_{m-1}.
&
\nonumber
\end{eqnarray}
Additionally, define  $\bar{a}_i:\mathbb{R}^2\to \mathbb{R}$  by setting
\begin{eqnarray}
& \bar{a}_i({\zeta}_1,{\zeta}_2)=\langle a_{ij}\delta_{j1}\efr_*\rangle {\zeta}_j,\ \ i,j=1,2. 
\label{DfBr-a-LnFn}
\end{eqnarray}
\begin{dfntn}\label{DfSnrSlwSst}
\textsf{\emph{By the  slow system of order $N=1,2,\ldots$ we mean the system of equations relative to the mean values  of corrections of order $N$, $\bar{p}_N=\bar{p}_N(x,t)$, $\bar{s}_N=\bar{s}_N(x,t)$, that reads
\begin{eqnarray}
&
\bar{p}_{N,t}+\nabla\cdot\left(\bar{p}_{N}\overline{\vc}+\chi\bar{p}\nabla \bar{s}_N-\mu{\bar{\mathcal{D}}}\nabla \bar{p}_{N}\right)
=\bar{a}_{1}(\bar{p}_{N},\bar{s}_{N}) +\bar{y}_{1,N},
&
\label{EqPrdSlwOrdr>2}\\
&
\bar{s}_{N,t}=\bar{a}_2(\bar{p}_{N},\bar{s}_{N})+\bar{y}_{2,N},
&
\label{EqPrySlwOrdr>2}
\end{eqnarray}
where the notations of $\bar{\mathcal{D}}$ and  $\bar{\vc}$ are for  the tensor and vector  fields that we had already defined in \eqref{DfDrftAndVscTnsrLdng},
\begin{eqnarray}
&
\bar{y}_{2,N}=\langle  {a}_2(\pch^\circ_{N},\widetilde{s}_{N}) +\sum\limits_{j=1..N-1}s_jg_{N-j}\rangle-\Delta \bar{s}_{N-1},
&
\nonumber\\
&
\bar{y}_{1,N}=\langle{a}_1(\pch^\circ_N,\widetilde{s}_{N})+ \sum\limits_{j=1..N-1}p_jf_{N-j}-\pch^\circ_{N,t} \rangle+
&
\nonumber\\
&
+\nabla\cdot\left\langle\widetilde{\nabla}_\efr\mathcal{L}^{-1}(\sum\limits_{j=0..N-1}p_jf_{N-1-j}-
\pch^\circ_{N,\tau} -p_{N-1,t}-\nabla\cdot\vq_{N-1})\right\rangle+
&
\nonumber\\
&
\nabla\cdot\left\langle(I-\mathcal{P})(\chi\sum\limits_{j=1}^N(\pch_j\widetilde{\nabla}\widetilde{s}_{N+1-j}
-\pch^\circ_j{\nabla}{s}_{N-j})-\mu\nabla_\efr \pch_N
)\right\rangle,
&
\nonumber
\end{eqnarray}
and equality \eqref{EqFlxCfcOrdr>0-1} with $k=N$ defines  vector $\vq_{N-1}$. 
}}
\end{dfntn} 
\textbf{Remark 8.} The  slow equations of order $N$  represent the expressions  linear  in the corrections of the same order. This is true not only regarding the projections considered as the unknowns,  $\bar{s}_N$ or $\bar{p}_N$, but also for  those seen as known. These are $\widetilde{s}_N$, $\pch_{N}$, $\pht_N^\circ$, which enters  the right hand sides, $\bar{y}_{N,1}$ and $\bar{y}_{N,2}$. The last ones  depend explicitly neither on the unknowns,  nor on $\pht^\bullet_N$, which is the  counterpart of $\bar{p}_N$.
  \begin{theorem}\label{PrpsSnrSlwSst} 
\textsf{\emph{Let  $N\ge 1$. Assume that we know the  $k-$th order corrections,  $s_k,p_k$, 
for  every $k=0..N-1$. Additionally, assume that system \eqref{EqPrdCfcOrdr>0}-\eqref{EqPryCfcOrdr>0} is compatible for $k=N$ and  determines  the projections of  $N-$th order corrections, $\pht^\circ_{N}$, $\pch_N$ and $\widetilde{s}_{N}$ uniquely, while projections  $\pht^\bullet_{N},\bar{s}_N$ remains undetermined and  the mean values, $\bar{p}_N$, $\bar{s}_N$ deliver a solution  to the  slow system of order $N$.   Then  system  \eqref{EqPrdCfcOrdr>0}-\eqref{EqPryCfcOrdr>0} is compatible for $k=N+1$ and determines projections of  the next order corrections, $\pht^\circ_{N+1}$, $\pch_{N+1}$ and $\widetilde{s}_{N+1}$ uniquely, while leaving projections  $\pht^\bullet_{N+1},\bar{s}_{N+1}$ undetermined. 
}}
\end{theorem}
\textbf{Remark 9.}  The formulation of Theorem~\ref{PrpsSnrSlwSst} presumes the independence of projection $\pht^\circ_{N}$, $\pch_N$ and $\widetilde{s}_{N}$ of the  complementing projections, $\pht^\bullet_{N},\bar{s}_N$. We consider the latter one  as one of the unknowns in the slow system, while the former one is  the image of  the other unknown upon  one-to-one mapping.   So,  the slow system's right hand sides, $\bar{y}_{N,1}$ and $\bar{y}_{N,2}$, depend on the unknowns neither explicitly nor implicitly. This observation makes  the slow system of order $N$ linear.  
\\ 
{\bf Proof of Theorem~\ref{PrpsSnrSlwSst}.}   For $k=N+1$, the   compatibility condition for equation   \eqref{EqPryCfcOrdr>0}  reads exactly as slow equation \eqref{EqPrySlwOrdr>2} of order $N$. 
Since this equation holds true by the assumption,  we  conclude that
 \begin{eqnarray*}&
 s_{N+1}=\bar{s}_{N+1}+\widetilde{s}_{N+1},\quad \widetilde{s}_{N+1}=\mathcal{H}^{-1}\left(\sum\limits_{j=1..N-1}s_jg_{N-j}-H \widetilde{s}_{N}-\Delta s_{N-1}\right),
 &\end{eqnarray*}
 where  $\mathcal{H}^{-1}:\ker M\to \ker M$. Hence, projection $\widetilde{s}_{N+1}$ is uniquely defined.  
 
 Equation   \eqref{EqPrdCfcOrdr>0}  also needs a  compatibility condition    to hold for $k=N+1$.  This one reads 
 \begin{eqnarray}&\label{EqCmptPrdOrd>0}
   \langle p_{N+1} \rangle^\xi_{\tau}+\left\langle\sum\limits_{j=0..N}p_jf_{N-j} - p_{N,t}-\nabla\cdot\vq_{N}\right\rangle^\xi=0
 &\end{eqnarray}
 We consider this as an equation for  determining  the spatial average of the next correction, $p_{N+1}$.  Further,  we eliminate vector $\vq_N$ by setting 
\begin{eqnarray}
&
\vq_N=\widetilde{\nabla}_{\efr}\mathcal{L}^{-1}\phi_{N}
+(I-\mathcal{P})\vu_N,\quad \text{where}\quad \vu_{N}={\hat{\vu}}_{N}+\check{\vu}_{N},
&
\label{EqQNDcmps}
\\
&
\phi_{N}=\sum\limits_{j=0..N-1}p_jf_{N-1-j} -p_{N-1,t}-\nabla\cdot\vq_{N-1}-{p}_{N,\tau},\quad
&
\nonumber\\
&
\hat{\vu}_{N}= \chi\sum\limits_{j=0}^N\pht_j\nabla\bar{s}_{N-j} -\mu (\nabla \pht_{N}-\pht_{N}\nabla\ln\langle\efr\rangle^\xi),
&
\nonumber\\
&
\check{\vu}_{N}=\chi\sum\limits_{j=1}^N\pch_j(\widetilde{\nabla}\widetilde{s}_{N+1-j}-{\nabla}{s}_{N-j})-\mu\nabla_\efr \pch_N.
&
\nonumber
\end{eqnarray}
where we have been using equation \eqref{EqFlxCfcOrdr>0-1} for $k=N$ to express vector $\vq_{N-1}$ and  
we have been taking into account  that 
$
p_0=\pht_0^\bullet, \pch_0=\pht^\circ=0.
$
The  action of  right inverse operator, $\mathcal{L}^{-1}$, involved on the right hand side in equality \eqref{EqQNDcmps} is well-defined, since $M^\xi\phi_N=0$, that is a corollary to the compatibility of system \eqref{EqPrdCfcOrdr>0}-\eqref{EqPryCfcOrdr>0}  for $k=N$, which, in turn, holds true  by assumption. At the same time, equality \eqref{EqQNDcmps}  is consistent with equation \eqref{EqPrdCfcOrdr>0} for $k=N$. Indeed, $\widetilde{\nabla}\cdot (I-\mathcal{P})=0$ by construction, hence,
\begin{eqnarray*}&
\widetilde{\nabla}\cdot\vq_N=\widetilde{\nabla}\cdot\widetilde{\nabla}_{\efr}\mathcal{L}^{-1}\phi_N=\phi_N=
{\sum\limits_{j=0..N-1}}p_jf_{N-1-j} -p_{N-1,t}-\nabla\cdot\vq_{N-1}-p_{N\tau}.
&\end{eqnarray*}
Now we aim at matching expression \eqref{EqQNDcmps} with equation \eqref{EqFlxCfcOrdr>0} for $k=N+1$ by choosing a suitable function ${p}_{N+1}$. It leads to the equation as follows
 \begin{eqnarray*}
&
-\mu \widetilde{\nabla}_\efr {p}_{N+1}+\chi{\sum\limits_{j=0}^{N}}\pht_j\widetilde{\nabla} \widetilde{s}_{N+1-j}  + \vu_N=\widetilde{\nabla}_{\efr}\mathcal{L}^{-1}\phi_{N}
+(I-\mathcal{P})\vu_N.
&
\nonumber\
 \end{eqnarray*}
The equivalent form of the last one reads
 \begin{eqnarray*}
&
\mu \widetilde{\nabla}_\efr {p}_{N+1}+\chi\efr_*\widetilde{\nabla}\sum\limits_{j=0}^{N}\bpx_j \widetilde{s}_{N+1-j} -\widetilde{\nabla}_{\efr}\mathcal{L}^{-1}\phi_{N}
+\mathcal{P}\vu_N=0.
&
\nonumber
 \end{eqnarray*}
From here, by a piece of argument similar to that involved in the proof of Theorem~\ref{PrpsLdngSlwdSst}, we isolate $p_{N+1}$ up to its projection on $\mathrm{range}\, P$, while determining $\pch_{N+1}$ uniquely. The result reads 
\begin{eqnarray*}
&
\mu  {\pch}_{N+1}=\chi Q\sum\limits_{j=0}^{N}\pht_j\widetilde{\nabla} \widetilde{s}_{N+1-j} +\mathcal{L}^{-1}(\phi_{N}
-\widetilde{\nabla}\cdot\vu_N).
&
\nonumber
 \end{eqnarray*}
The projection on $\mathrm{range} P$ is free to choose still, and we aim at using this to take care of the compatibility of equation \eqref{EqCmptPrdOrd>0} (where vector $\vq_N$ is specified by formula \eqref{EqQNDcmps}). Its compatibility condition reads 
\beq 
\left\langle \sum\limits_{j=0..N}p_jf_{N-j} - p_{N,t}-\nabla\cdot\vq_{N}\right\rangle=0,
\eeq
where formula \eqref{EqQNDcmps}  gives expression for  vector $\vq_{N}$. From here, the direct calculations with taking into account definition  \eqref{DfMtrxM} eventually  bring us at equation~\eqref{EqPrdSlwOrdr>2}, which holds by assumption. From this point, by repeating  a suitable piece of argument from the proof of  Theorem~\ref{PrpsLdngSlwdSst}, we arrive at the single-valued definition of component   $\pht^\circ_{N+1}$ by formula
\beq
\pht^\circ_{N+1}=P\pr_\tau^{-1} \left\langle\sum\limits_{j=0..N}p_jf_{N-j} - p_{N,t}-\nabla\cdot\vq_{N}\right\rangle^\xi,
\eeq 
where operator  $\pr^{-1}$ is the right inverse  that sends $\mathrm{range}\,(I-M^\tau)$ to itself,  while vector $\vq_N$ is specified by equality \eqref{EqQNDcmps}. 

The component $\pht^\bullet_{N+1}$ remains indefinite.  Thus, we have determined the corrections of order $N+1$, $p_{N+1}$ and $s_{N+1}$, up to projections $\pht^\bullet_{N+1}$ and $\bar{s}_{N+1}$, and this completes the proof. 
 \\$\square$\\
\textbf{Remark 10.} A direct inspection of the last proof  shows that components  $\pch_{N+1}$, $\pht^\circ_{N+1}$ and $\widetilde{s}_{N+1}$ are $Op_N-$expressions linear in the corrections of order $N$, $p_{N}$ and ${s}_{N}$. Hence, developing the successive  asymptotic approximations in question is a linear process except for the leading approximation, as the leading slow system is non-linear.
\section{Imposed Stability and Instability}\label{ScStb}
\noindent 
Following to  Kapitza's theory of the upside-down pendulum \cite{LndLfs}, we see the imposed stability or instability  as the stability or instability  of the special solutions called quasi-equilibria.  Every quasi-equilibrium  allows for identifying with a suitable equilibria, and we can address the stability of the latter one applying the common techniques. 
\subsection{Quasi-equilibria}\label{SscQsEq}
\noindent
The quasi-equilibria patterns (if any) are  the short waves, which  propagate with no slow modulation at least in the leading approximation -- that is, by formulae \eqref{P0} and \eqref{S0}, $s_0\equiv s_e$, $p_0\equiv p_e\efr_*$, where $p_e=\const$ and $s_e=\const$. Therefore, it is of sense to identify the quasi-equilibria with  the  equilibria of the leading slow system -- that is, with its special solutions, such that 
  \begin{eqnarray}&\label{QsEqlbr}
 \bar{s}=\bar{s}_e\equiv\const,\quad\bar{p}=\bar{p}_e \equiv \const.
&\end{eqnarray}
 Feasibility of the  homogeneous equilibria  is natural  for a  system possessing the total translational invariance. The leading  slow systems possess this invariance  provided that the signal undergoes no slow modulation -- that is, $h=h(\xi,\tau)$. 
The  quasi-equilibria densities, $\bar{p}_e$ and $\bar{s}_e$,  must obey the following equations
\begin{eqnarray} &\label{BrPeSe}
  \bar{p}_e\bar{f}(\bar{p}_e,\bar{s}_e)=0,\quad\bar{s}_e\bar{g}(\bar{p}_e,\bar{s}_e)=0,
& \end{eqnarray}  
where the notations of  $\bar{f}$ and $\bar{g}$ are for functions defined in \eqref{EqBrP} and \eqref{EqBrS}. Generally, they differ from the original kinetic terms,  $f$ and $g$.  An interesting feature is that $\bar{f}=f$ and $\bar{g}=g$ provided that the original kinetic terms are  linear in $p$ (briefly, $p-$linear).  Then the quasi-equilibria coincide with the usual equilibria. In particular, Lotka-Volterra's  or   Holling's II and III kinetics (see \cite{TtnTrFn} for discussion)  possess this feature.
\subsection{Linear stability analysis}\label{SscLnStbAnls}
\noindent
In the sequel,  we assume that the external signal does not undergo a slow modulation -- that is, $h=h(\xi, \tau)$ by default.

To examine the  stability of a quasi-equilibrium, we perform so-called linear stability analysis  for the corresponding equilibrium -- that is,  we linearize  the leading slow system nearby  it,
and then look for a special solution  that reads 
\begin{equation}\label{EgnMds}
      (\hat{p},\hat{q},\hat{s})\exp(i k x+\lambda t),\ \hat{p},\hat{q},\hat{s},\lambda \in \mathbb{C},\quad k\in \mathbb{R}^n,
\end{equation}
 It is common to name solutions \eqref{EgnMds} as the normal modes (of the corresponding equilibrium) and say that such a mode is \textsl{stable~(unstable, neutral)} if the real part of spectral parameter, $\lambda$,  is negative (positive, equal to zero). 
  It is also common to say that an equilibrium is stable~(unstable, neutral) provided every its normal mode is stable (there exists an unstable (neutral) mode). 
  
  For the normal modes, the so-called spectral stability problem arises for which the spectral parameter $\lambda$ is the eigenvalue, and the corresponding normal mode is the eigenmode. Moreover,  the separation of the variables reduces the spectral stability problem to the algebraic eigenvalue one for every concrete wave vector, $k$.
  
  Let the signal be issued with some characteristic temporal frequency, $c$, with some characteristic amplitude $a_0$ --  that is,
  \beq
  h=a_0 h_0(\xi,c\tau)
  \eeq
   and function $h_0$ is 1-periodic in variable $c\tau$ with unit norm  in space $\mathrm{L}_2(\mathbb{T}^{n+1})$, say.  Then  
  \beq
   \efr=\exp(a_0\varkappa h_0/\mu),
   \eeq 
 so we consider number  $a=a_0\varkappa/\mu>0$ as   an effective amplitude. Further, set
  \beear
  &\hat\mu=\mu k\cdot\bar{\mathcal{D}} k,\ \quad \hat\chi=p_e\chi k\cdot\bar{\mathcal{D}} k,\quad  \hat{c}=\bar{\vc}\cdot k&
  \label{DfHtPrmt}
  \eear
  where the notations of $\bar{\mathcal{D}}$ and  $\bar{\vc}$ are for  the tensor and vector  fields that we had already defined in \eqref{DfDrftAndVscTnsrLdng}. It worths noticing, that 
  \beq
  \hat{c}=c\bar{\vc}_1\cdot k, 
\eeq  
  where vector $\bar{\vc}_1$ is the normalized drift velocity  evaluated in the same way as vector $\vc$, but for  $c=1$. 
  
 It turns out that all the statements on page 12 in \cite{AM3} remains valid in the main respects upon replacing the parameters defined there by those defined above. 
 In more details, for  a quasi-equilibrium with densities $\bar{p}_e$, $\bar{s}_e$ consider an eigenmode with wave vector $k$. Consider a triad of real numbers
 \begin{eqnarray}
\Delta_0=1,  &\Delta_2=\hat{\mu}+\hat{\delta}-{\bar{a}^\circ}_{{1,1}}-{\bar{a}^\circ}_{{2,2}},\quad \Delta_4= \left({\hat{c}}^{2}+\Delta_2^2\right)(\hat{\delta} -{\bar{a}^\circ}_{{2
,2}})\left(\hat{\mu}-{\bar{a}^\circ}_{{1,1}} \right)-{\bar{a}^\circ}_{{2,1}}\Delta_2^2({\bar{a}^\circ}_{1,2}+\hat{\chi})&
\label{HnklMnrsGnrl}
 \end{eqnarray}
 where the notations of ${\bar{a}^\circ}_{{i,j}}$, $i,j=1,2$, are for the coefficients defined by formulae \eqref{EqDf-a-LnFn}-\eqref{DfBr-a-LnFn} for ${p}_0=\bar{p}_e\efr_*$, $\bar{s}=\bar{s}_e$. The number of the sign changes in this triad relative to its current ordering is equal to the number of unstable eigenmodes linked to the equilibria and wave vector chosen. In particular, the quasi-equilibrium remains linearly  stable  as long as 
\begin{equation}\label{RstrBrAij}
{\bar{a}^\circ}_{11}\le 0,\ {\bar{a}^\circ}_{22}\le 0,\ {\bar{a}^\circ}_{21}\le 0,\ {\bar{a}^\circ}_{12}\ge 0
\end{equation}
where at least 3 inequality are strict. 

Suppose  we deal with a $p-$linear kinetics. Given the comment put at the end of Subsection~\ref{SscQsEq}, we see that ${\bar{a}^\circ}_{{i,j}}={{a}^\circ}_{{i,j}}$, $i,j=1,2$, where the notations of ${{a}^\circ}_{{i,j}}$ are for the coefficients defined by formulae \eqref{EqDf-a-LnFn} for $p_0=p_e$ and $s_0=s_e$, and numbers  $p_e$ and $s_e$ are the values of the equilibrium densities -- that is,
\beq
{p}_e\bar{f}(\bar{p}_e,\bar{s}_e)=0,\quad\bar{s}_e\bar{g}(\bar{p}_e,\bar{s}_e)=0.
\eeq
(currently, the equilibria and quasi-equilibria coincide). Thus,  every equilibrium in every system \eqref{EqPrdTrns}-\eqref{EqPryDnst} with $p-$linear kinetics is stable either with or without external signal provided that 
\begin{equation}\label{RstrAij}
{{a}^\circ}_{11}\le 0,\ {{a}^\circ}_{22}\le 0,\ {{a}^\circ}_{21}\le 0,\ {{a}^\circ}_{12}\ge 0
\end{equation}
Although there is no  diffusion  in the prey equation of the leading slow system,  let's allow for it, but set it weak, for better understanding its effect. Let the notation of $\delta$ be  for the small diffusivity coefficient.
\begin{lemma}\label{LmmHnklChn} 
\textsf{\emph{Let the parametric point 
$(\hat{\mu},\hat{\delta},\hat{\chi},{\bar{a}}^\circ_{11},{\bar{a}}^\circ_{12},{\bar{a}}^\circ_{21},{\bar{a}}^\circ_{22})$ 
satisfy  inequalities
\begin{eqnarray}
&{\bar{a}}^\circ_{11}{\bar{a}}^\circ_{22}-{\bar{a}}^\circ_{21}{\bar{a}}^\circ_{12}>0,\ {\bar{a}}^\circ_{11}+{\bar{a}}^\circ_{22}<0,\  {\bar{a}}^\circ_{11}{\bar{a}}^\circ_{22}<0,&\ 
\label{InstLmmCnd1}\\
&(\hat{\delta} -{\bar{a}}^\circ_{{2
,2}})\left(\hat{\mu}-{\bar{a}}^\circ_{{1,1}} \right)<0,\ (\hat{\delta} -{\bar{a}}^\circ_{{2
,2}})\left(\hat{\mu}-{\bar{a}}^\circ_{{1,1}} \right)-{\bar{a}}^\circ_{{2,1}}({\bar{a}}^\circ_{1,2}+\hat{\chi})>0,&
\label{InstLmmCnd2}\\
&\hat{\mu}\ge 0,\hat{\delta}\ge 0,\hat{\chi}\ge 0.&
\nonumber
\end{eqnarray}
 Then, for every such a parametric point, there exists a threshold value, 
\begin{eqnarray}&\label{ThrsholdCSqrd}
  \hat{c}_*^2=-\Delta_2^2\frac{(\hat{\delta} -{\bar{a}}^\circ_{{2
,2}})\left(\hat{\mu}-{\bar{a}}^\circ_{{1,1}} \right)-{\bar{a}}^\circ_{{2,1}}({\bar{a}}^\circ_{1,2}+\hat{\chi})}{(\hat{\delta} -{\bar{a}}^\circ_{2,2})(\hat{\mu}-{\bar{a}}^\circ_{1,1} )},
&\end{eqnarray} 
   such that the chain of signs in triad \eqref{HnklMnrsGnrl}  reads
\begin{eqnarray*}
  &+,+,+ \ \text{for}\ \hat{c}^2<\hat{c}^2_* &
  \\
  &+,+,-  \ \text{for }\ \hat{c}^2>\hat{c}^2_*.&
\end{eqnarray*}
}}
\end{lemma}
Further, set  
\beq
\alpha\byd |k|,\quad k=\alpha\bar{k},\quad |\bar{k}|=1.
\eeq
 \begin{theorem}\label{ThInst}
 \textsf{\emph{
 Given a system  \eqref{EqPrdTrns}-\eqref{EqPryDnst} in multiple dimensions with a short-wave external signal, $h$, issued on  characteristic frequency $c$ with effective  amplitude $a$.   Let this signal undergo no slow modulation and determines nonzero vector $\vc$.  Let coefficients $a^\circ_{ij}$ be computed 
 by formulae \eqref{EqDf-a-LnFn} for $p_0=p_e$ and $s_0=s_e$.
 Assume the following
\begin{description}
  \item[(i)] the kinetics given allows for a non-degenerated equilibrium, $p=p_e,s=s_e$;
  \item[(ii)] all the conditions of Lemma~\ref{LmmHnklChn} remain  true if we put  $a^\circ_{ij}$ in place of $\bar{a}^\circ_{ij}$.
\end{description}
Then there exist $a_*>0$ and the branch of quasi-equilibria $\bar{p}_e=\bar{p}_e(a)$, $\bar{s}_e=\bar{s}_e(a)$,  which is smooth on $[0,a_*)$,  such that $\bar{p}_e(0)=p_e$,  $\bar{s}_e(0)={s}_e$, and for every  $a\in [0,a_*) $ and for every unit vector $\bar{k}:\ \vc_1\cdot\bar{k}\neq 0$ there exists $k_*\in(0,\infty]$ (we do not exclude $k=+\infty$) such that for every  wave vector $k=\alpha\bar{k}: 0<\alpha< k_*$ there exists a threshold value $c_*=c_*(\bar{k},\alpha,a)\ge 0$ such that the normal mode linked to quasi-equilibrium $\bar{p}_e(a)$, $\bar{s}_e(a)$  and to wave vector $k$  is stable  for $c^2<c^2_*$ and unstable for $c^2>c^2_*$.  
 }}
\end{theorem}
 \textbf{ Remark 11. } The threshold value of the characteristic frequency, $c_*$ goes  to infinity when $\alpha=|k|\to 0$ uniformly in directors,  $\bar{k}$. This excludes the longwave instability. In particular, the homogeneous mode remains stable.
\section{Discussion And Conclusions}\label{ScCncl}
\noindent
Thus, we have developed the full asymptotic expansion for the  solutions to the predator-prey system with Keller-Siegel's responding to the  prey density  and to an  external signal in multiple dimensions with no special assumptions on the external signal or solution. This result extends the asymptotic of short traveling waves in one dimension \cite{AM3} to the general shortwave signal  and multiple dimensions. In one spatial dimension, we arrive  at the same leading slow system as that derived in \cite{AM3} (see subsection \eqref{Ssc1Dmn}). By that, we have freed the conclusions regarding the imposed stability made therein  of the special assumptions on the  form of the solutions.

Although the leading term of an asymptotic  gives us a summary of the most substantial information on the  solution,  the higher order terms  can play an important part   in  justifying the asymptotic approximation.  Regarding the asymptotic presented above, such a justification seems to be  quite feasible because only the leading slow system is non-linear, while the construction of the subsequent approximations is a purely linear process. It involves smoothing in the fast variables, and losing the derivatives in the slow ones. The latter feature is easy to see from the equations for the  mean prey density in the high order slow systems. These  are simply linear ODEs in the time variable, but with right hand sides involving the derivatives  of the preceding approximations in the spatial variables.  So,   calcualting  only a finite order  asymptotic approximation is feasible assuming a finite regularity for the problem data.

Further, we have studied the imposed stability due to a general external signal in multiple dimensions. We have  proved that a signal is unable to  impose the instability to the homogeneous equilibria if the kinetics is linear in the predators density ($p-$linear) and inequalities \eqref{RstrAij} hold true. (Remind, by the way, that the signal is  unable to modify $p-$linear kinetics, and the quasi-equilibria coincide with the equilibria, therefore.) In mathematical ecology, the  $p-$linear  kinetics  often reads  
\begin{equation*}
  f= f(s),\ g=g_0(s) -g_1(s)p/s,\quad \text{where}\ f_s>0,\ g_{0s}<0,\ g_1>0\ \forall\,s>0
\end{equation*}
For instance,  letting $g_1=s$, $g_0=1-s$, $f(s)=\gamma g_1-\beta$, we get Lotka-Volterra's kinetics if $g_1=s$, and we get Holling's II-III if $g_1=s^n/(1+bs^n)$, n=1,2,\ldots.  Additionally, let's assume  that
\begin{eqnarray}
&
g_{0s}-\frac{sg_0}{g_{1}}\left(\frac{g_{1}}{s}\right)_s< 0\quad \forall\,s>0:\ g_0(s)>0.
&
\label{PrCptRtnIncrs}
\end{eqnarray}
Then such a kinetics obeys inequalities \eqref{RstrAij}. Straightforward examination  shows that inequalities \eqref{PrCptRtnIncrs} hold true for Lotka-Volterra's kinetics, and for Holling III too.  For Holling-II kinetics, if we put it into the dimensionless form  as follows $g_1=\alpha s/(1+\alpha s)$, inequalities \eqref{PrCptRtnIncrs} hold  true  provided that $\alpha<1$.

The Arditi-Ginzburg ratio-dependent kinetics (\cite{ArG}, and see  discussions in \cite{TtnTrFn}) turns out to be capable of behaving  with more diversity. This kinetics reads (upon a suitable scaling) 
\beq
  f=\frac{\gamma\,s}{s+p}-\beta;\quad g=1-s-\frac{rp}{s+p},\ r>0,
\eeq
where number $r$ plays the part of the parameter. In a suitable parametric domain it  allows for the coexistence of both species in equilibrium. There is a parametric sub-domain where this equilibrium   gets unstable with no signal. Regarding this feature we allude  to  works \cite{LeeHlnLws},\cite{BnPtrRtDpnd}. Since we are interested in the instabilities imposed by the external signal,  let's restrict the parameters in within a domain  where the coexistence equilibrium is stable provided that the signal is off. However,   this time  the signal is able to modify the kinetics -- that is, $f\neq \bar{f}$, and $g\neq \bar{g}$. Since the study of this effect is beyond the score of this article,  we restrict ourselves within  the signals of small effective amplitude, $a$. It turned out that even such a weak signal can  impose instability provided that its characteristic frequency, $c$, is great enough. Earlier, article \cite{AM3} have reported this, but for one dimension and with all other restrictions mentioned above. The speed at which the signal propagates has been playing the part of the  characteristic frequency.

Overall, the present study acknowledges one more time that  the external signal is likely not capable of creating  the instability domain in the parametric space from nothing but it can substantially  widen the one that is non-empty with no signal, as it has been supposed in \cite {AM3} relying on the one dimensional analysis.  The novel feature, however, is the substantial dependence on the   direction of the wave vector, $k$, and, especially, on the angle between the wave vector and the drift velocity.  It is very clear, for example, from   the criterion for instability delivered by Lemma~\ref{LmmHnklChn}.  Among the normal modes  with equal magnitudes of $|k|$,  more dangerous (in the sense of getting unstable) are those linked to the wave vectors parallel to the drift velocity, $\vc$ (provided it is not zero, of course). For example,  the signal that propagates as a traveling wave in the direction specified by vector $\theta$ produces the drift towards the same direction (see Subsection~\ref{Ssc1Dmn}).

Prediction of the drift  in the full generality is far beyond the scope of the present study. This also concerns the diffusivity tensor resulting from the homogenization of the general shortwave signal.  The explicit examples are still feasible only for a special class of signals (see Subsection~\ref{Ssc1Dmn}), so the   prediction and control of these fields represent an open mathematical problem doubtlessly  relevant to the practical issues.

{\bf Acknowledgments.} The authors are thankful to  Southern Federal University for the opportunity to do this research. By the same reason, A. Morgulis is grateful to the South Mathematical Institute of the Vladikavkaz Scientific Center of RAS.

\end{document}